
\input amstex
\documentstyle{amsppt}
\refstyle {A}
\widestnumber\key{AAAAA}
\magnification=1200
\parskip 6pt
\pagewidth{5.25in}
\pageheight{7.56in}
\def\leq{\le} 
\def\geq{\ge} 
\def\CC{{\Bbb C}}
\def\ZZ{{\Bbb Z}}
\def\CD{{\Cal D}}
\def\CE{{\Cal E}}
\def\CF{{\Cal F}}
\def\CJ{{\Cal J}}
\def\CV{{\Cal V}}
\def\CZ{{\Cal Z}}

\def\eps{\varepsilon}
\def\O{{\Cal O}}

\def\phi{\varphi}
\def\\{\par\noindent}

\def\mult{\text{mult}}
\def\dim{\text{dim}}
\def\ind{\text{ind}}
\def\lra {{\longrightarrow}}

\def \I {\Cal I}


\topmatter
\title Bounds for Seshadri constants  \endtitle
\author Oliver K\"{u}chle$^*$  and Andreas Steffens$^{**}$ \endauthor
\address{ Oliver K\"{u}chle \hfill {} \linebreak
\hglue 13.6pt Mathematisches Institut  \hfill {} \linebreak
\hglue 13.6pt Universt\"at Bayreuth\ \hfill {} \linebreak
\hglue 13.6pt D--95440 Bayreuth \hfill {} \linebreak
\hglue 13.6pt Deutschland}  \endaddress
\email  kuechle\@btm8x1.mat.uni-bayreuth.de \endemail
\address{ Andreas Steffens \hfill {} \linebreak
\hglue 13.6pt Department of Mathematics \hfill {} \linebreak
\hglue 13.6pt UCLA \hfill {} \linebreak
\hglue 13.6pt Los Angeles, Ca 90024 \hfill {} \linebreak
\hglue 13.6pt USA}  \endaddress
\email  steffens\@math.ucla.edu \endemail
\thanks
\vglue1pt
$^*$Supported by Max--Planck--Gesellschaft and Deutsche
Forschungsgemeinschaft\endthanks
\thanks $^{**}$Supported by Deutsche Forschungsgemeinschaft \endthanks
\rightheadtext{}
\leftheadtext{}
\endtopmatter

\document
\head  Introduction \endhead
In this paper we present an alternative approach to the boundedness
of Seshadri constants of nef and big line bundles
at a general point of a complex--projective variety.

Seshadri constants $\eps(L,x)$, which have been introduced
by Demailly \cite{De92}, measure the local positivity of
a nef line bundle $L$ at a point $x\in X$
of  a complex--projective variety $X$,
and can be defined as
$$\eps(L,x):=\inf_{C\ni x}\left\{\frac{L\cdot C}{\text{mult}_x(C)}\right\},$$
where the infimum is taken over all reduced irreducible curves
$C\subset X$ passing through $x$ (cf. also \S1 below, \cite{De92} or
\cite{EKL} for further
characterizations and properties of Seshadri constants).

Over the last years there has been quite some activity in
studying Seshadri constants,
starting with the somehow surprising result
by Ein and Lazarsfeld (cf. \cite{EL}) that the Seshadri constant of
an ample line bundle on a smooth surface is bounded below by
$1$ for all except perhaps countably many points.

On the other hand examples by Miranda (cf. \cite{EKL, 1.5})
show that for any
integral $n\ge 2$ and real $\delta >0$
there is a smooth $n-$dimensional variety $X$, an ample
line bundle $L$ on $X$ and a point $x\in X$ with
$\eps(L,x) < \delta$; in other words, there does
not exist a universal lower bound for Seshadri constants
valid for all $X$ and ample $L$ at every point $x\in X$.

Then it was proven by Ein-K\"uchle-Lazarsfeld \cite{EKL}
that, for a nef and big line bundle $L$ on an $n-$dimensional
projective variety, the Seshadri constant at very general
points (i.e. outside a countable union of proper subvarieties)
is bounded below by $\frac{1}{n}$, and  that this implies
the existence of a lower bound at general points depending
only on $n$.

Finally we want to mention the recent papers by
Nakamaye \cite{Na} and Steffens \cite{St}
dealing with the problem of
"maximality" of Seshadri constants on abelian varieties, as well
as variants due to K\"uchle \cite{K\"u} and Paoletti \cite{Pa}
concerning Seshadri constants along several points and
higher dimensional subvarieties respectively.

Here we first reprove the existence of lower bounds for Seshadri constants
of nef and big line bundles at general points of
a projective variety.
Although the bound obtained does, in general, not improve
the one in \cite{EKL},
the method at hand may give better results in certain cases,
since our bound can be expressed more flexibly in terms
of the degrees of subvarieties
with respect to the line bundle in question:
\proclaim{Theorem 3.4}
Let $L$ be a nef line bundle on an $n$-dimensional
irreducible projective variety $X$ and
$\eps >0$ be a real number.
Let $\alpha_1,\dots,\alpha_{n}$ be positive rational numbers
and $x\in X$ a general point. Put $\gamma =1+\sum_{i=1}^{n-1}\alpha_i$.
Suppose that any $d-$dimensional
$(1\le d\le n)$ subvariety $V\subseteq X$ containing a very general
point $y\in X$ satisfies
$$\text{deg}_LV=
L^d.V \ge \frac{\eps^d\cdot
\gamma^n}
{\alpha_d^{n-d}}.$$
Then $\eps(L,x)\ge \eps.$
\endproclaim
One special feature of our proof is that there is no
mention of "Seshadri--exceptional" curves;
the objects of consideration are families of divisors
with high multiplicity at prescribed points.
The method itself is based upon ideas of Demailly's paper
\cite{De93} as explained and translated in an algebro--geometric language
by Ein--Lazarsfeld--Nakamaye in \cite{ELN}.
The basic idea is, very roughly,
to start with an effective divisor $E$
in a linear series $|kL|$ ($k\gg0$) with large multiplicity
at $x$, and to consider the locus $V$ of points where
the singularities of $E$ are "concentrated" in a certain way.
The possibility that $V$ is zero--dimensional
imposes constraints on
the local positivity of $L$ at $x$ in a sense. Otherwise one uses
variational techniques to give a lower bound on the
degree of $V$.

Instead of, as in \cite{De93} and \cite{ELN}, making a positivity
assumption on the tangent bundle of the manifold in question
to be able to "differentiate", we apply the strategy of
differentiation in parameter directions in the spirit of
\cite{EKL}.
This, however, forces us to consider only very general points
instead of arbitrary ones.
The result of this method is the following Theorem, which
might be of interest also in other contexts:
\proclaim{Theorem 2.1}
Let $X$ be
a smooth n-dimensional projective variety, $L$  a nef and big line
bundle on $X$ and $\alpha >0$ a rational number such that
$L^n>\alpha^n$. Let
$$0=\beta_1<\beta_2<\dots <\beta_n<\beta_{n+1}=\alpha $$
be  any  sequence
of rational numbers and $x\in X$ a very  general point. Then
 either
\roster
\item
"$(a)$" there exist $k\gg 0$ and a divisor $E\in|kL|$ having an
isolated singularity with multiplicity at least
$k(\beta_{n+1}-\beta_{n})$ at $x$, or
\medskip
\item"$(b)$"
there exists a proper subvariety $V\subset X$ through $x$ of codimension
$c \leq n-1$  such that
$$\deg_LV=L^{n-c}\cdot V\leq\frac{1}{(\beta_{c+1}-\beta_c)^c}
\left(1-{\root n \of {\left(1-\frac{\alpha^n}{L^n}\right)^c}}\right)L^n
<\frac{\alpha^n}{(\beta_{c+1}-\beta_c)^c}.$$
\endroster
\endproclaim
To pass from
Theorem 2.1 to actually bounding the
Seshadri constant we use a rescaling--trick (cf. 3.3)
in combination with the well known characterization of Seshadri
constants via the generation of $s-$jets by certain
adjoint linear systems (cf. 1.5).

After a first version of this paper was written, we realized
that the rescaling argument mentioned above
can also be applied to the results of \cite{ELN},
leading to bounds for Seshadri constants valid
at {\sl arbitrary} points; these bounds, however,
depend on the line bundle $L$ and the manifold $X$,
or, rather, its tangent bundle $T_X$:
\proclaim{Corollary 4.4}
Let $X$ be a smooth $n$--dimensional projective variety,
$x\in X$ any point,
$A$ an ample line bundle and $\delta\ge 0$ a real number
such that $T_X(\delta A)$ ist nef.
Then
$$\eps(A,x)\ge \text{min}\left\{\frac{1}{(n-1)^{n-1}(2n-1)^n},
\frac{1}{\delta}\right\}.$$
\endproclaim
This gives in particular bounds valid at arbitrary points for
the Seshadri constants of
the canonical line bundle $K_X=(\bigwedge^nT_X)^*$ or
its inverse if these are ample, or for any ample line bundle
in case $K_X$ is trivial (cf. 4.5).

Corollary 4.4 is in accordance with and should be compared to bounds
following from recent work of Angehrn--Siu \cite{AS} on the
basepoint--freeness of adjoint linear series (cf. 4.6).

The paper is organized as follows. After fixing notations
and establishing a general setup in \S0 we recall some
basic facts about Seshadri constants and collect some
auxiliary statements in \S1. Then, in \S2, we
prove the main technical result, Theorem 2.1.
Finally we give the applications to bounding Seshadri constants
at general points in \S3, and at arbitrary points in \S4.

\subhead
Acknowledgements
\endsubhead
We are grateful to Rob Lazarsfeld for not only having
drawn our attention to Seshadri constants but also
for sharing his insights in this topic
with us. In particular, he suggested to combine the ideas
of \cite{ELN} with the methods of \cite{EKL} towards one of us.
The reader will notice that our treatment relies on
ideas and prior works due to Demailly \cite{De93},
Ein--K\"uchle--Lazarsfeld \cite{EKL},
Fujita \cite{Fu94}, and especially
Ein--Lazarsfeld--Nakamaye \cite{ELN}
to whom we are indebted.

\head \S0. Notations and the General Setup \endhead

\subhead
(0.1)
\endsubhead
Throughout this paper we will work over the field $\CC$ of
complex numbers. Given a variety $Y$ a {\sl general} point $y\in Y$
is a point in the complement of some proper
subvariety of $Y$, and a {\sl very general} point is a point in
the complement of some countable union of proper subvarieties.

\subhead
(0.2)
\endsubhead
Given a smooth variety $Y$, an integer $m\ge 0$ and
a subvariety $W\subset Y$
we denote by $\I_W^{(m)}$ the symbolic
power sheaf of all functions vanishing to order at least $m$ along $W$.
Then $\I_W^{(1)}=\I_W$ is the ideal sheaf of
$W$, and $\I_Z^{(m)}=\I_Z^m$ for a smooth subvariety $Z\subset Y$.

\subhead
(0.3)
\endsubhead
Let $M$ be a line bundle on a smooth variety $Y$.
Given a divisor $E\in |kM|$ we call the normalized multiplicity
$$\ind_y(E)=\frac{\mult_y(E)}{k}$$
the {\it index} of $E$ at a point $y\in Y$. \\
We say that a divisor $E\in |kM|$ has an {\it almost isolated singularity of
index $\geq\alpha$} at a point $y\in Y$, if
$\ind_y(E)\geq \alpha$ and
there is an open neighborhood $U\subset Y$ of $y$ such that
$$ \ind_x(E) < 1 \text{ for all } x\in U\setminus \{y\}.$$

\subhead
(0.4)
\endsubhead
For a line bundle $L$ and a coherent sheaf of ideals $\CJ$ on $Y$
we denote by $|L\otimes\CJ|$ the linear subsystem of the complete
linear system $|L|$ corresponding to sections in $L\otimes \CJ$.
Given such a system $|L\otimes \CJ|\ne\emptyset$ on $Y$, by
the {\it base locus} Bs$|L\otimes \CJ|$ we mean the support
of the intersection of all members of $|L\otimes\CJ|$.

\subhead
(0.5)
\endsubhead
We will be concerned with the following setup:
Let be $X$ a smooth irreducible
$n$--dimensional projective variety,
$T$ a smooth irreducible affine variety and
$g: T\lra X$ a quasi-finite dominant morphism with graph $\Gamma$.
Let $pr_X$ and $pr_T$ denote the projections from $Y=X\times T$
onto its factors.
Note that the projections then map $\Gamma$ dominantly to $X$ resp. $T$.
Given a Zariski-closed subset
(or subscheme) $Z \subset X \times T$, we consider
the fibre $Z_t$ of $pr_T$
as a subset (or subscheme) of $X$. Similarly, $Z_x \subset T$ is the
fibre of $Z$ over $x \in X$.
Given a sheaf $\CF$ on $X\times T$ we write $\CF_t$ for the induced
sheaf on $X$.

Examples of this situation are provided
by Zariski--open subsets $T\subset X$.

\head  \S1. Preliminaries \endhead

In this section we collect some preliminary results which will be needed
in the sequel.
We start with some remarks concerning multiplicity loci in a
family.

Let $E$ be an effective divisor on a smooth variety $Y$.
Then the function $y \to \mult_y(E)$ is Zariski upper-semicontinuous on $Y$.
For any given irreducible subvariety $Z\subset Y$ we refer with $\mult_Z(E)$
to the value of $\mult_z(E)$ at a general point $z\in Z$.
The following lemma allows one to make
fibrewise calculations of multiplicities.

\proclaim{Lemma 1.1} (cf. \cite{EKL, 2.1})
Let $X$ and $T$ be smooth irreducible varieties, and suppose that
$Z\subset X\times T$ is an irreducible subvariety which dominates $T$.
Let $\CE \subset X\times T$ be an effective divisor.

Then for a general point $t\in T$, and any irreducible component
$W_t\subset \CZ_t$ of the fiber $\CZ_t$, we have
$$ \mult_{W_t}(\CE_t)= \mult_{\CZ}(\CE).\qed $$
\endproclaim

The next elementary lemma gives a way to detect
irreducible components of base loci.

\proclaim{Lemma 1.2}
Let $m$ be a positive integer,
$M$ a line bundle on a smooth variety $Y$
and $V\subseteq W\subseteq Y$
subvarieties such that $V$ is irreducible.
Suppose that
\roster
\item
$\I_W\otimes M$ is generated by global sections, and
\smallskip
\item
$\I_W\subseteq \I_V^{(m)}$.
\endroster

\noindent
Then $V$ is an irreducible component
of Bs$|\I_V^{(m)}\otimes M|$.
\endproclaim
\demo{Proof}
Consider the inclusions
$$V\subseteq \text{Bs}|\I_V^{(m)}\otimes M|\subseteq\text{Bs}|\I_W
\otimes M|\subseteq W,\tag{*}$$
where the first inclusion is clear, the second
follows from (2), and the third from the
fact that Bs$|\I_W\otimes M|=W$ because
$\I_W\otimes M$ is globally generated according to (1).

Now let $Z\subseteq  \text{Bs}|\I_V^{(m)}\otimes M|$ be any irreducible
component containing $V$. Then  $(*)$ and the assumption that
$V\subseteq W$ is an irreducible component imply
$V=Z$,  hence the claim follows.
\qed
\enddemo

Now we show that being an irreducible subvariety is well-behaved in
families:
\proclaim{Lemma 1.3}
Let $f:Y \lra Z$ be a morphism between irreducible
varieties, and $V, W\subseteq Y$ be subvarieties such
that $V\subseteq W$ as an irreducible subvariety and $V$
maps dominant to $Z$.

Then over general $z\in Z$ every irreducible component $U_z$ of $V_z$
has dimension dim$V -\text{dim}Z$ and is an irreducible
component of $W_z$.
\endproclaim
\demo{Proof} For the first part we refer to \cite{Ha, II.Ex 3.22}.
The second assertion comes down to an easy dimension count as follows:
Write $$W=V\cup V'$$
with a subvariety $V'\subseteq Y$ not containing $V$, i.e. dim$V >
\text{dim}(V\cap V')$. If $V\cap V'$ does not map dominantly to $Z$,
then $V_z$ and $V'_z$ do not meet in a general fibre and the claim
follows.

Otherwise we obtain over general $z\in Z$:
$$\align \text{dim}V_z &= \dim V-\dim Y\cr
		       &>\dim(V\cap V')-\dim Y\cr
		       &=\dim(V\cap V')_z,\cr\endalign$$
and therefore $\dim U_z >\dim (V\cap V')_z\ge\dim (U_z\cap V'_z)$,
where $U_z\subseteq V_z$ is any irreducible component.
Considering the decomposition $V_z=U_z\cup U'_z$,
where $U_z\not\subseteq U'_z$, we conclude
$\dim U_z > \dim (U_z\cap(U'_z\cup V'_z))$,
which shows that $U_z$ is an irreducible component
of $W_z=U_z\cup(V'_z\cup U'_z)$.
\hfil\qed
\enddemo

For the reader's convenience we recall some well--known facts
concerning Seshadri constants. The next Lemma deals with the
Seshadri constant at general versus the one at very general points
(cf. \cite{EKL, 1.4}).
\proclaim{Lemma 1.4}
Let $X$ be a smooth projective variety and $L$ a nef and big
line bundle on $X$. Suppose $\eps(L,y)=\eps$ for a point $y\in X$.
Then for any real $\delta > 0$ there exists exists a Zariski--open
subset $U(\delta)\subseteq X$ such that
$$\eps(L,x)\ge \eps-\delta\text{ for all } x\in U(\delta).\hfil\qed$$
\endproclaim

Finally we recall the relations between almost isolated
singularities, generation of higher jets, and Seshadri constants
(cf. \cite{ELN, 1.1}, \cite{EKL, 1.3}).

\proclaim{Theorem 1.5}  Let $X$ be a smooth projective variety
of dimension $n$ and $L$ a nef and big line bundle on $X$.

\parskip 5pt
\noindent {\rm (1.5.1).}
Suppose there exists a divisor $E \in |kL|$
having an almost isolated singularity of index $\ge n+s$ at $x\in X$.
Then $$H^1(X, \O_X(K_X + L) \otimes \I^{s+1}_x)  = 0.$$
In particular, $|K_X + L|$ generates $s$-jets at $x$,
i.e. the evaluation map
$$H^0(X, \O_X(K_X + L)) \lra  H^0(X, \O_X(K_X +
L) \otimes \O_X / \I_x^{s+1})$$
is surjective.

\parskip 5pt
\noindent {\rm (1.5.2).}
Let $\eps(L,x)$ be the Seshadri constant of $L$ at $x$.
If $${\displaystyle r > {s\over \eps(L,x)} + {n\over
\eps(L,x)}},$$  then $|K_X+rL|$ generates $s$-jets at $x\in X$. The same
statement holds if
 ${\displaystyle r = {s +n\over \eps(L,x)}}$ and $L^n>\eps(L,x)^n$.

\parskip5pt
\noindent {\rm (1.5.3).}
Conversely, suppose there is a real number $\eps >0$ plus a real constant $c$
such that $|K_X+rL|$ generates $s$-jets at $x\in X$ for all $s\gg 0$
whenever
$${\displaystyle  r > {s\over \eps}+c}.$$  Then $
\eps(L,x)\ge \eps.$

\noindent
In other words: If $\eps(L,x)\ < \alpha$, then for all
$s_0 >0$ and all real $c$ there exist an $s\ge s_0$ and an $r >
\frac{s}{\alpha} +c$ such that $|K_X+rL|$ does not generate $s-$jets at $x$.
\qed
\endproclaim

\head \S2.  Relative  moving part estimates \endhead
The aim of the section is to proof the following Theorem.

\proclaim{Theorem 2.1}
Let $X$ be
a smooth n-dimensional projective variety, $L$  a nef and big line
bundle on $X$ and $\alpha >0$ a rational number such that
$L^n>\alpha^n$. Let
$$0=\beta_1<\beta_2<\dots <\beta_n<\beta_{n+1}=\alpha $$
be  any  sequence
of rational numbers and $x\in X$ a very  general point. Then
 either
\roster
\item
"$(a)$" there exist $k\gg 0$ and a divisor $E\in|kL|$ having an
isolated singularity of index
$\ind_x(E)\geq \beta_{n+1}-\beta_{n}$ at $x$; or
\medskip
\item"$(b)$"
there exists a proper subvariety $V\subset X$ through $x$ of codimension
$c \leq n-1$  such that
$$\deg_LV=L^{n-c}\cdot V\leq\frac{1}{(\beta_{c+1}-\beta_c)^c}
\left(1-{\root n \of {\left(1-\frac{\alpha^n}{L^n}\right)^c}}\right)L^n.$$
\endroster
\endproclaim

\subhead (2.2) Families of divisors \endsubhead
Pick an arbitrary point $y\in X$ and a smooth affine
neighbourhood $T\subset X$ of $y$ in $X$. Then the embedding
$g: T\hookrightarrow X$ satisfies the properties of (0.5).
Note that (very) general points of $T$ correspond to (very) general
points of $X$. We will use the notations introduced in (0.5) henceforth.

Argueing as in \cite{EKL, (3.8)}, for $k\gg 0$ with $\alpha k\in \ZZ$
we obtain divisors $\CE_k\in |pr_X^*(kL)|$
in $X\times T$ satisfying
$$\ind_\Gamma(\CE_k) > \alpha.$$
The argument is, in brief, that for any
$x\in X$  using  Riemann-Roch and a parameter count one finds
a divisor $E\in  |kL|$ with  $\mult_x(E) > \alpha k$.

Hence the torsion free $\O_T$-module
$$R=pr_{T*}\left(pr_X^*(kL)\otimes \I_\Gamma^{(\alpha k)}\right)$$

has positive rank,
and is globally generated since $T$ is affine.
Therefore a non-zero section $\Phi$  of $R$  gives via the evaluation map
$pr_T^*R \to pr_X^*(kL)\otimes \I_\Gamma^{(\alpha k)}$ the desired
divisor $\CE_k$.

\subhead (2.3) The  multiplicity schemes $\CZ_\sigma(\CE)$
\endsubhead
For $k\ge 0$ with $\alpha k\in \ZZ$ put
$$A_k = |\I_\Gamma^{(\alpha k)}\otimes pr_X^*(kL) |.$$
Then by (2.2) $A_k$ is non--empty  for sufficiently large $k$.
For nonzero $\CE_k\in A_k$ and  any rational number $\sigma > 0$
we define
$$\CZ_\sigma(\CE_k)=
\{\ y=(x,t)\in \CE_k\ | \ \ind_y(\CE_k) \geq \sigma\  \}.$$
Note that $\CZ_\sigma(\CE_k)$ is a Zariski-closed
subset of $X\times T$. Its  natural
scheme structure is given locally by the vanishing of all partial derivatives
of order $< k\sigma $ of a local equation for $\CE_k$.
We will be only interested in $\CZ_\sigma(\CE_k)$ as an algebraic
set, and for a general
choice of $\CE_k$.

Recall the following generalized version of Bertini's Theorem,
due to Koll\'ar (cf. \cite{Ko}):
\proclaim{Theorem 2.3.0}
Let $Y$ be a smooth variety and $|B_1,\dots,B_k|$ be
a linear system on $Y$. Then a general member $B\in |B_1,\dots,B_k|$
satisfies
$$\mult_yB\le 1+ \inf_i\left\{\mult_y B_i\right\}$$
at every point $y\in Y$.
\qed\endproclaim
The following lemma, which is an analog
of \cite{ELN,(3.8)}, says that the multiplicity loci
$\CZ_\sigma(\CE_k)$
are independent of $k$  and the choice of
a general $\CE_k\in A_k$ as
soon  as $k$ is sufficiently large; here by general
we mean general in the sense of Theorem 2.3.0.
\proclaim{Lemma 2.3.1}
For fixed $\sigma$ there is a positive integer $k_0$
such  that $\CZ_\sigma(\CE_{k_1})=\CZ_\sigma(\CE_{k_2})$ for all
$k_1, k_2
\geq k_0$.
\endproclaim
\demo{Proof (cf. \cite{ELN})}
 As $\sigma$ and $\alpha$ are fixed, we
will for simplicity write $\CZ(k)$ for $\CZ_\sigma(\CE_k)$.
Choose an integer $m\geq 2$ such that
$A_k\ne\emptyset$ for $k\ge m$.
Fixing an integer $a \ge m$ we then claim that there exists a
positive integer
$k(a)$ such that
$$ \CZ(c) \subseteq \CZ(a) \ \text{ whenever } \ c \ge k(a). \tag{*}$$
To prove the claim,  suppose that $y \not \in \CZ(a)$ so that there
exists $\eta > 0$ satisfying
$$\mult_y(\CE_a) \le a \sigma - \eta.$$
Note that since the index is a discrete invariant,
$\eta$ is bounded below independently of $y$; in fact if $m \sigma \in \ZZ$
then
$\eta \geq 1/m$. Suppose $b\ge m$ is an integer relatively prime to
$a$.  Then any integer
$c \ge ab$ can be expressed as
$$ c = \alpha a + \beta b, \ \ \alpha, \beta \in \ZZ \ \text{and}
\ 0 \le \beta \le a.
$$
Consider the divisor
$\CE_c^\prime = \alpha \CE_a + \beta \CE_b \in A_c.$ Then
$$\aligned
\mult_y(\CE_c)\leq 1+\mult_y(\CE_c^\prime) &
= 1+\alpha \cdot \mult_y(\CE_a) + \beta \cdot
\mult_y(\CE_b) \\
                &\le 1+\alpha a \sigma - \alpha \eta + \beta \cdot
                    \mult_y(\CE_b) \\
                & = c \left ( \sigma -
                \frac{\eta\left(1 - \frac{\beta b}{c}\right)}{a}
     + \frac{1+\beta \cdot \mult_y(\CE_b)}{c} \right ) ,
\endaligned
$$ where the first inequality is a consequence
of Theorem 2.3.0.
Since $\eta, \beta,$ and $b$ are bounded independently of $c$, it
follows that
$\mult_y(\CE_c) < c \sigma$ for $c \gg 0$. Hence $y \not \in \CZ(c)$ for
all sufficiently large $c$ as claimed.

If $\CZ(c) = \CZ(a)$ for all $c \gg 0$ then we are finished.  If not, then
by $(*)$ there exists  $a^\prime > 0$ such that $\CZ(a^\prime)
\subsetneq \CZ(a)$. The argument can then be repeated with
$a^\prime$ instead of $a$.  This process cannot go on indefinitely and
this establishes the Lemma.
\hfil
\qed
\enddemo
The next Lemma is an adaptation of
\cite{ELN, (1.5),(1.6)} to our situation.
We present a sketch of its proof for the reader's convenience.
\proclaim{Gap-Lemma 2.3.2}
Let $\CE\subset X\times T$  be a
family of effective divisors on $X$ with $\ind_\Gamma(\CE) > \alpha$ along
the graph $\Gamma\subset X\times T$ of $g:T\to X$.
Let
$$0=\beta_1< \beta_2 < \dots < \beta_n < \beta_{n+1} =\alpha$$
be any sequence  of rational numbers.
Define $$\CZ_0=X\times T \text{ and }
\CZ_j=\CZ_{\beta_j}(\CE)=\{\ y\in\CE\ | \ \ind_y(\CE)\geq \beta_j\ \}
\text{ for } 1\leq j\leq n+1.$$
Then there exists an index $c,\ 1\le c\le n$, and an irreducible subvariety
$\CV\subset
X\times T$
such  that:
\roster
\item
$codim(\CV)=c$,

\item
 $\Gamma\subset\CV$, and

\item
 $\CV$ is an irreducible component of both
 $\CZ_c$ and $\CZ_{c+1}$.
\endroster
\endproclaim
This means that the index of $\CE$ "jumps" by at least $\beta_{c+1}-\beta_c$
along
$\CV$, i.e.
$\ind_y(\CE)\geq \beta_{c+1}$ for every $y\in\CV$ and there is
an open set $U\subset X\times T$ meeting $\CV$ such that
$\ind_v(\CE) < \beta_c$ for every $v\in U\setminus \CV$.

\demo{Sketch of Proof} The sets $\CZ_i$ lie in a chain
$$\Gamma\subseteq \CZ_{n+1} \subseteq \dots
\subseteq \CZ_1=\CE\subsetneq \CZ_0=X\times T.$$
Starting with $\CZ_{n+1}$ and working up in dimension, we can choose
irreducible components $\CV_j$ of $\CZ_j$ containing $\Gamma$ such that
$\CV_{j+1} \subseteq \CV_j$. So we arrive at a chain of irreducible
varieties
$$
\Gamma\ \subseteq
\ \CV_{n+1} \ \subseteq \ \CV_n \ \subseteq \ldots \ \subseteq \CV_1 \
\subsetneq \CV_0 \ = \ X\times T,$$
and since $X\times T$ is irreducible of dimension
$2n=\text{dim}(\Gamma)+n$,
at least two consecutive links in the chain must coincide,
say $\CV_c = \CV_{c+1}$, and we take $\CV = \CV_c$.
Using elementary combinatorial arguments one shows that also the
condition $codim(\CV)=c$ can be achieved. For details
we refer to the proof of Lemma (1.6) in \cite{ELN}.
\hfil
\qed
\enddemo

\subhead
(2.4)
\endsubhead
{}From now on we fix the set of rational numbers
$$0=\beta_1< \beta_2 < \dots < \beta_n < \beta_{n+1} =\alpha.$$

Then by Lemma 2.3.1 there exists an integer $k_0$ such that
the multiplicity loci $\CZ_{\beta_i}(\CE_k),$ $1\le i\le n+1,$ are
independent of $k$ as soon as $k\ge  k_0$
and $\CE_k\in A_k$ are general.
Therefore also the multiplicity "jumping" loci $\CV$
obtained by the Gap-Lemma 2.3.2 can be chosen
independently of $\CE=\CE_k$ and $k$ up to the above restrictions.
Fix such a $\CV$ and put $\beta=\beta_{c+1}-\beta_c$.

\proclaim{Proposition 2.4.1} For all sufficiently divisible
$k\gg 0$
the jumping locus
$\CV$ is an irreducible component of the base locus of the linear system
$$| \I_{\CV}^{(\beta k)} \otimes pr_X^*(kL)|.$$
\endproclaim
Here by sufficiently divisible we mean that
$\beta_i k \in \ZZ$ for all $i=1,\dots,n+1$.

We start by recalling some general facts concerning the
differentiation of sections of line bundles $pr_X^*(kL)$
in parameter directions and its connection to certain
multiplicity loci (cf. also \cite{EKL, \S 2} and \cite{ELN, \S 2}).

Let $\CD^\ell_{X\times T}(pr_X^*(kL))$ be the sheaf of differential
operators of
order $\leq \ell$ on $pr_X^*(kL)$ and let $\CD^\ell_T$ be the be the sheaf
of differential operators  of order $\leq \ell$  on $T$.
Since there is a canonical inclusion of vector bundles
$$pr_T^*(\CD^\ell_T) \hookrightarrow \CD^\ell_{X\times T}(pr_X^*(kL)),$$
the sections of $\CD^\ell_T$ act naturally on the space of sections of
$pr_X^*(kL)$.
A section $\psi \in \Gamma(X\times T, pr_X^*(kL))$ determines a
homomorphism of sheaves
$${\frak d}_\ell(\psi) : pr_T^*(\CD^\ell_T) \to pr_X^*(kL).$$
Locally, represent $\psi$ by a function $f$, then ${\frak d}_\ell(\psi)$ is
just
taking a differential operator $D$ to the function $D(f)$. Since
$pr_X^*(kL)$ is a
line bundle there exists a sheaf of ideals $\I_{\Sigma_\ell(\psi)}$ such that
$$Im({\frak d}_\ell(\psi) : pr_T^*(\CD^\ell_T) \to pr_X^*(kL))=
\I_{\Sigma_\ell(\psi)}\otimes pr_X^*(kL).$$
Let $\psi$ be a defining section for a divisor $\CE\in |pr_X^*(kL)|$.
Then we claim that
$$\Sigma_\ell(\psi)=
\{\ (x,t)\in X\times T \ |\ \mult_t(\CE_x) > \ell\ \}.\tag{2.4.2}$$
Indeed, the scheme structure on the right hand side is given
locally by the vanishing of all partial derivatives in $T$ direction
of order $\leq\ell$ of a local equation for $\CE$.

Note also that the $\I_{\Sigma_\ell(\psi)}\otimes pr_X^*(kL)$ are as
quotients of the globally generated sheaf $pr_T^*\left(\CD^\ell_T\right)$
generated by global sections.

\demo{Proof of the Proposition}
The plan is to apply Lemma 1.2.
By assumption
$$\alpha k,\ \
p =\beta_c k \text{ and }
q = \beta k \text{ are integers}.$$
Let $\CE\in A_k$ be a divisor determining $\CV$ and $\psi$ a section
defining $\CE$. For integers $\ell$ put $\Sigma_\ell=\Sigma_\ell(\psi)$.
We claim that
$$\I_{\Sigma_{p-1}} \subset \I_\CV^{(q)}.\tag{2.4.3}$$
To prove this let
$f$ be a local equation for $\CE$ over some open set $U$.
Then $\I_{\Sigma_{p-1}}$ is locally
generated by all functions
$$\{\ D(f)\ |\  D\in (pr_T^*\CD^{p-1}_T)(U)\ \}.$$
On the other hand we have that $\tilde{D}(f)\in \I_\CV$ for every
$\tilde{D}\in\CD^{p+q-1}_{X\times T}(pr_X^*(kL))(U)$,
since $\CV$ is an irreducible component of
$\CZ_{\frac{p+q}{k}}(\CE)$. And in particular $\tilde{D}(f)\in \I_\CV$ for
every
$\tilde{D}\in (pr_T^*\CD^{p+q-1}_T)(U)$. Hence for all $D\in
pr_T^*\CD^{p-1}_T$
the function $D(f)$ vanishes to the order $\geq q$ on $\CV$ which
shows (2.4.3).

Since $\I_{\Sigma_{p-1}}\otimes pr_X^*(kL)$ is globally generated,
application of Lemma 1.2 to our situation will give the desired result
once we show that $\CV$ is also an irreducible component of $\Sigma_{p-1}$.
This is the content of Lemma 2.4.4 below.
\hfil
\qed
\enddemo

\proclaim{Lemma 2.4.4} Let $\sigma k$ be a positive integer,
$\CE\in |pr_X^*(kL)|$ be an effective
divisor on $X\times T$ and $V\subset \CZ_\sigma(\CE)$
an irreducible component
dominating $X$. Then $V$ is also an irreducible component of
$\Sigma_{k\sigma -1}(\CE)$.
\endproclaim
\demo{Proof}
By definition $\CZ_\sigma(\CE)\subset \Sigma_{k\sigma -1}(\CE)$.
Let $W\subset\Sigma_{k\sigma-1}(\CE)$ be an irreducible component
containing $V$. If we can show that
$W\subset  \CZ_\sigma(\CE)$, then we are done, because that implies
$V=W$. Lemma 1.1 shows $\ind_W(\CE)=\ind_{W_x'}(\CE)$ for general
$x\in X$ and any irreducible component $W_x'$ of $W_x$.
Hence the assertion follows from
$\ind_{W_x}(\CE_x) \geq \sigma$, where we used (2.4.2).
\hfil\qed\enddemo

\proclaim{Corollary 2.4.5} For all sufficiently  divisible $k\gg 0$,
we have
$$| \I_\Gamma^{(\alpha k)}\otimes pr_X^*(kL)| \subset
| \I_{\CV}^{(\beta k)} \otimes pr_X^*(kL)|.$$
\endproclaim
\demo{Proof}
By the above any sufficiently general $\CE\in A_k$ determines
the same $\CV$, in particular such an $\CE$ satisfies
$\ind_\CV(\CE)\ge\beta_c$,
and this implies $\CE\in |\I_{\Sigma_{p-1}}\otimes pr_X^*(kL)|$
by (2.4.2), where again we assume that $\alpha k,\ p=\beta_c k$ and
$q=\beta k$ are integral.
Then the claim follows from (2.4.3).\hfil
\qed
\enddemo
\proclaim{Proposition 2.4.6}
For all sufficiently divisible $k \gg 0$ and
very general $t\in T$ the following hold:
\roster
\item
There exists an irreducible subvariety $V\subseteq X$
of codimension $c$ containing $g(t)$ which is an irreducible
component of the base locus of the linear system
$$|\CJ_k|:=\left|\left(\I_{\CV}^{(\beta k)}\right)_t\otimes kL\right|$$
on $X$ such that mult$_VD\ge k\beta$ for all $D\in |\CJ_k|$.
In particular, if $c=n$, i.e. $\CV=\Gamma$, then
by Bertini's Theorem there exists
a divisor $D\in |kL|$ having an isolated
singularity of index $\ge \beta$ at the point $x=g(t)$.
\medskip
\item
$\dim H^0(X, \I_{g(t)}^{\alpha k}\otimes kL )\leq
\dim H^0(X,\CJ_k).$
\endroster
\endproclaim
\demo{Proof}
To begin with we study the situation for a fixed $k$.
First we note that after possibly shrinking $T$
we can assume that the coherent sheaf
$\CF=\I_{\CV}^{(\beta k)}\otimes pr_X^*(kL)$ is flat over $T$.
In fact $pr_T: X\times T\lra T$ is projective and $T$ is affine
and integral, hence the assertion follows from consideration
of the Hilbertpolynomials of the $\CF_t$ (cf. also \cite{Ha, III.9.9}):
These do not depend on $t$ for $t$ in an open dense subset of $T$.

After possibly shrinking $T$ more it follows from semicontinuity
that there is a natural isomorphism
$$H^0(X\times T,\CF)\otimes k(t)\simeq H^0(X,\CF_t).\tag{$*$}$$
(cf. \cite{Ha, III.12.9}).
In other words, taking  global sections of $\CF$ commutes with restricting to
fibres over general $t\in T$, and therefore
$(\text{Bs}|\CF|)_t=\text{Bs}|\CJ_k|$ for such $t$.

Now we can prove (1).  By Proposition 2.4.1 we have
$\Gamma\subseteq \CV\subseteq \text{Bs}|\CF|$ with
$\CV$ an irreducible component of Bs$|\CF|$, hence Lemma 1.3
shows that any irreducible component $V$ of
$\CV_t$  is an $(n-c)-$dimensional
irreducible component of $(\text{Bs}|\CF|)_t=\text{Bs}|\CJ_k|$.
It remains to show $\mult_VD\ge k\beta$ for all $D\in |\CJ_k|=|\CF_t|$.
But this follows from $(*)$ and Lemma 1.2.

Assertion (2) follows in the same way from Corollary 2.4.5
and the fact that $$
\left(\I_{\Gamma}^{(\alpha k)}\right)_t=
\left(\I_{\Gamma}^{\alpha k}\right)_t=\I_{g(t)}^{\alpha k}.$$
To complete the proof of the Proposition we only have to
remark that, since $\CV$ does not depend on $k$,
the above arguments work simultaneously for
all divisible $k\gg 0$ if we replace the general $t\in T$
by a very general $t\in T$.
\hfil\qed
\enddemo

\subhead (2.5) Bounding the degree of irreducible components of base loci
\endsubhead
In this subsection we complete the proof of Theorem 2.1 by
bounding the degree of the irreducible component $V$ from 2.4.6
using a strategy essentially due to Fujita (cf.  \cite{Fu82}, \cite{Fu94}).
Alternatively one could carry out an approach
via graded linear series
as in \cite{ELN} leading to slightly weaker  bounds.

Let $k\gg 0$ be sufficiently large and divisible, and fix a very
general $t\in T$.
Let $\CJ_k=(\I_{\CV}^{(\beta k)})_t\otimes kL$
and $V\subset X$ be as in (2.4.6); recall that $V$ depends
on $t$ but {\sl not} on $k$, and that
$V$ is an $(n-c)$--dimensional irreducible component
of Bs$|\CJ_k|$.
We  may assume that $\dim V >0$,
since otherwise the assertion
of Theorem 2.1 follows from (2.4.6.1).
\subhead
(2.5.1)
\endsubhead
Resolving the base locus of $\Lambda:=|\CJ_k|$  we can find a sequence
$$X^\prime =X_s\to\dots\to X_r\to X_{r-1}\to\dots\to X_1\to X_0=X$$
of birational morphisms $\tau_i:X_i\to  X_{i-1}$
together with linear systems $\Lambda_i$ on $X_i$
such that
\roster
\item
$\Lambda_0=\Lambda$.
\smallskip
\item
$\tau_i: X_i\to X_{i-1}$ is the blow-up of a smooth
subvariety $C_i$ of $X_{i-1}$.
\smallskip
\item
$\tau_i^*\Lambda_{i-1}=\Lambda_i + m_iE_i$ for some
nonnegative integers $m_i$,
where $E_i$ is the exceptional divisor on $X_i$ lying over $C_i$ and
$E_i\nsubseteq \text{Bs}\Lambda_i$.
\smallskip
\item
Bs$\Lambda_s=\emptyset$.
\endroster
Let $\tau=\tau_s\circ\dots\circ\tau_1$ be the composition,
$E^*_i$ resp. $E^\prime_i$ be the total resp. proper
transforms of $E_i$ on $X^\prime$, and
$Y_i=\tau(E^*_i)=\tau(E^\prime_i)$.
The $Y_i$  coincide with the image of $C_i$ in $X$.

Let $F_i$ denote  the pull back to $X^\prime$ of the
general member of
the linear system $\Lambda_i$ on $X_i$.
Finally, let $H$ be a general member of $\Lambda_s$,
$F=\tau^*(kL)$, and $E$ the fixed part of $\tau^*\Lambda$, so
that $H=F-E$  where  $E=\sum_{i=1}^{s}m_iE^*_i$.

By assumption there exists an index $r$ with $Y_r=V$, and
$$m_r\ge k\beta\tag{2.5.2}$$
since $\mult_V(D)\geq k\beta$ for all $D\in|\CJ_k|$.
We also may and will assume that the resolution $\tau$
is chosen in such a way that
$\dim(Y_i) < \dim(V) =n-c$ for all $i<r$.

\proclaim{Lemma 2.5.3}
$F^{n-c}.E.H^{c-1}\ge k^{n-c}m_r^c \deg_LV$.
\endproclaim
\demo{Proof}
First of all note that, since $F$ and $H$ are nef and
$E-m_rE^*_r$ is an effective divisor, we have
$$F^{n-c}.E.H^{c-1}
\geq m_r F^{n-c}.E_r^*.H^{c-1}.\tag{*}$$
Next we use that
$$F^{n-c}.E^*_r.H^{c-1}=F^{n-c}.E_r^*.F_r^{c-1},\tag{**}$$
which follows from the standard Segre class computations
because of
$$\text{dim}\left(\text{supp}\left(\tau_*\left( E^*_r.F^{c-1}_r -
E^*_r.H^{c-1}\right)\right)\right) < n-c.$$
Finally we compute $F^{n-c}.E_r^*.F_r^{c-1}$.
Since dim$(Y_i) < n-c$ for all $i<r$, we have
$$F^{n-c}.F_r^c=F^{n-c}.(F-\sum_{i\leq
r}m_iE^*_i)^c=F^{n-c}.(F-m_rE^*_r)^c,$$
hence $F^{n-c}.F_r^c= F^n - m_r^c\deg_{kL}(V)$ by the birationality of
the morphism $C_r\to Y_r=V$.
A similar argument shows $F^{n-c+1}.F_r^{c-1}=F^n$,
and therefore
$$\aligned F^{n-c}.E_r^*.F_r^{c-1}&=\frac{1}{m_r}\left(F^{n-c}F_r^{c-1}\right).
\left(F-F_r-\sum_{i<r}m_iE_i^*\right)\cr
&=\frac{1}{m_r}\left(F^{n-c+1}.F_r^{c-1}-F^{n-c}.F_r^c\right)\cr
&=m_r^{c-1}k^{n-c}\text{deg}_L V.\endaligned $$
Combining this with $(*)$ and $(**)$ proves the Lemma.\qed
\enddemo

\proclaim{Lemma 2.5.4} For any $\eps > 0$
there exist a sufficiently large and divisible
$k$ and a resolution $\tau: X'\lra X$ of the rational
map given by $|\CJ_k|$ satisfying the
properties in (2.5.1) and
$H^n=(F-E)^n\geq k^n(L^n-\alpha^n -\eps)$.
\endproclaim
\demo{Proof}
The proof follows closely the proof of the Theorem in \cite{Fu94}.
Therefore we only give an outline and indicate the necessary
modifications.
For varying $k$ consider $|\CJ_k|$ and
denote by $(X'_k,H_k)$ the pair
$(X^\prime,H)$ obtained as in (2.5.1).
We will derive a contradiction assuming that
$H_k^n < k^n(L^n - \alpha^n -\eps)$ on $X'_k$ for all
large and divisible $k$.

Letting $\eps$ grow if necessary, we may assume that
$H_\ell^n \geq \ell^n\left(L^n - \alpha^n - \eps- \frac{\eps}{(2n)!}\right)$
on $X'_\ell$ for one fixed large and divisible $\ell$.
Now, for any integer $s>0$, we claim that
$$h^0(X,\CJ_{s\ell}) \leq h^0(X'_\ell,sH_\ell)+\frac{n(s\ell)^n\eps}{(2n)!},$$
which is proven exactly as in \cite{Fu94} by
using the lower bound on $H^n_\ell$ and considering an
appropriate resolution of $\Lambda=|\CJ_{s\ell}|$.
{}From (2.4.6.2) and asymptotic Riemann--Roch we then obtain
$$\aligned
\frac{(s\ell)^n}{n!}\left(L^n-\alpha^n\right) + \phi(s)
&\leq h^0(X, I_x^{\alpha s\ell}\otimes s\ell L )\cr
&\leq h^0(X,\CJ_{s\ell)})\cr
&\leq h^0(X'_\ell, sH_\ell)+\frac{n(s\ell)^n\eps}{(2n)!}\cr
&\leq \frac{s^n}{n!}\ell^n\left(L^n-\alpha^n-\eps \right)
+\frac{n(s\ell)^n\eps}{(2n)!} +\psi(s),\cr
\endaligned$$
where $\phi$ and $\psi$ are functions with
$\lim_{s\to\infty}\frac{\phi(s)}{s^n}= \lim_{s\to\infty}
\frac{\psi(s)}{s^n}=0$. This gives the desired contradiction.
\qed
\enddemo

Before  stating the main result of this section which will
complete the proof of Theorem  2.1 we need to recall some well known
facts (cf. \cite{De93, 5.2}, \cite{Fu82, 1.2}).

\proclaim{Lemma 2.5.5}
Let $F$, $H$ be nef divisors on an $n$-dimensional smooth projective variety
$Y$.
Then:
\roster
\item
$F^d.H^{n-d} \geq {\root n \of {(F^n)^d}}{\root n \of {(H^n)^{n-d}}}$ holds
for all $0\leq d \leq n$.
\smallskip
\item
If $E= F-H$ is effective, then
$F^a.H^{n-a}\geq F^b.H^{n-b}$ for any $a\geq b$.
\qed\endroster\endproclaim

\proclaim{Proposition 2.5.6}
With the above notations the degree of $V$ satisfies
$$\deg_LV\leq
\frac{1}{\beta^c}
\left(1-{\root n \of {\left(1-\frac{\alpha^n}{L^n}\right)^c}}\right)L^n.$$
\endproclaim

\demo{Proof}
By Lemma 2.5.3 we have
$$\deg_LV \leq \frac{1}{k^{n-c}m_r^c}F^{n-c}.E.H^{c-1}.$$
Note that $F^{n-c}.E.H^{c-1}=\left(F^{n-c+1}.H^{c-1} - F^{n-c}.H^c\right)$.
So if
 we bound the first term using (2.5.5.2) and
bound the second term using (2.5.5.1), we find
$$\aligned
\deg_LV & \leq\frac{1}{k^{n-c}m_r^c}\left(F^{n-c+1}.H^{c-1} -
F^{n-c}.H^c\right)\cr
 & \leq\frac{1}{k^{n-c}m_r^c}
\left(F^n- (F^n)^{\frac{n-c}{n}}(H^n)^{\frac{c}{n}}\right) \cr
 & \leq\frac{k^n}{k^{n-c}m_r^c}\left(L^n-
(L^n)^{\frac{n-c}{n}}(L^n-\alpha^n)^{\frac{c}{n}}\right) \cr
 &\leq\frac{1}{\beta^c}
\left(1-{\root n \of {\left(1-\frac{\alpha^n}{L^n}\right)^c}}\right)L^n,
\endaligned$$
where the last steps are Lemma 2.5.4 plus the fact that deg$_LV$ is
integral, and (2.5.2).
\hfill
\qed
\enddemo


\head
\S3. Applications
\endhead

\proclaim{Theorem 3.1} Let $X$ be an $n$-dimensional smooth
projective variety, $x\in X$ be a very general
point and $L$ be a nef and big line bundle.
Let $r$ and $s$ be positive integers, $\gamma > 1$ a rational number
satisfying
$$(rL)^n > (\gamma(n+s))^n,$$ and $\alpha_1,\dots,\alpha_{n-1}$ be positive
rational numbers with $\sum_{i=1}^{n-1} \alpha_i=\gamma -1$.

\noindent
Then either

\noindent
$(a)\;\;\;\;|K_X+rL|$ generates $s$-jets at $x$, or

\noindent
$(b)\;\;\;$ there is a proper
subvariety $V\subset X$ of positive dimension $d$
containing $x$ with
$$\text{deg}_LV=
L^d.V\le\frac{(rL)^n}{ \alpha_d^{n-d}(n+s)^{n-d}r^d}\left(1-{\root n \of
{\left(1-\frac{\gamma^n(n+s)^n}{(rL)^n}\right)^{n-d}}}\right).$$
\endproclaim
\demo{Proof}
Put $$L'=rL,$$
$$\alpha=\beta_{n+1}=\gamma(n+s),$$
$$\beta_n=(\gamma -1)(n+s), \text{ and downward recursively}$$
$$\beta_i=\beta_{i+1}-\alpha_{n-i}(n+s)\;\ \text{for}\;\ i=n-1,\dots,1.$$
Let $x\in X$ be a very general point and apply Theorem 2.1.
Then either there exists $k\gg0$ and a divisor $E\in |kL'|$ having
an almost
isolated singularity of index $\ge n+s$ at $x$, in which case by (1.5.1)
the linear series $|K_X+L'|=|K_X+rL|$ generates $s-$jets at $x$, or
there exists a subvariety $V\subset X$ with the wanted properties.
\qed
\enddemo
\proclaim{Corollary 3.2}
With the assumptions of  Theorem 3.1 there exists $\kappa >0$ such that
either

\noindent
$(a)\;\;\;|K_X+rL|$ generates $s$-jets at $x$, or

\noindent
$(b)\;\;\;$ there is a proper subvariety $V\subset X$
of positive dimension $d$
containing $x$ with
$$\text{deg}_LV=
L^d.V < \left(\frac{n+s}{r}\right)^d \frac{\gamma^n}{\alpha_d^{n-d}}
\left(\frac{1}{1-\kappa}\right).$$
Given constants $c_1 > c_2 >0$, such $\kappa$ can be chosen independently
of $\frac{n+s}{r}, \gamma$ and $\alpha_d$  if the condition
$c_1\ge\frac{n+s}{r}\gamma\ge c_2$ holds.
\endproclaim
\demo{Proof}
This follows directly from Theorem 3.1 using the estimate
$$g\left(\frac{n+s}{r}\gamma;\kappa\right)
:=\left(1-\frac{(n+s)^n\gamma^n}{(rL)^n}
\left(\frac{1}{1-\kappa}\right)\right)^n
- \left(1-\frac{(n+s)^n\gamma^n}{(rL)^n}\right)^{n-d} < 0$$
for integral $1\le d < n$, fixed $\frac{n+s}{r}\gamma$,
and sufficiently small $\kappa$.
The independence of $\frac{n+s}{r}\gamma$ is a consequence of the fact
that $g$ satisfies  a Lipschitz condition with respect to
$\frac{n+s}{r}\gamma$.
\qed
\enddemo
\remark{Remark 3.3}
Up to the assumption on the positivity of $T_X$
and the genericity of $x\in X$ Theorem 3.1
looks similar to \cite{ELN, Theorem 4.1}. Note however that
in the estimate of the degree deg$_LV$ we have $(n+s)^d$
compared to $(n+s)^n$ in \cite{ELN}, which turns out to be crucial
when bounding the Seshadri constant.
This improvement is achieved by "rescaling" the intervall $[0,\alpha]$
from 2.1.
\endremark
\proclaim{Theorem 3.4}
Let $L$ be a nef line bundle on an $n$-dimensional
irreducible projective variety $X$ and
$\eps >0$ be a real number.
Let $\alpha_1,\dots,\alpha_{n}$ be positive rational numbers
and $x\in X$ be a general point.
Put $\gamma =1+\sum_{i=1}^{n-1}\alpha_i$.
Suppose that any $d-$dimensional
$(1\le d\le n)$ subvariety $V\subseteq X$ containing a very general
point $y\in X$ satisfies
$$\text{deg}_LV=
L^d.V \ge \frac{\eps^d\cdot
\gamma^n}
{\alpha_d^{n-d}}.$$
Then $\eps(L,x)\ge \eps.$
\endproclaim
\demo{Proof}
First we note that, since we are only considering general
points, there is no loss of generality in supposing that
$X$ is smooth (cf. \cite{EKL, 3.2} for the precise argument).
Suppose that $\eps(L,x) < \eps$ at general points $x\in X$.
Then by Lemma 1.4 for all $\delta >0$ we have $\eps(L,x) <\eps +\delta$
at any point $y\in X$,
and accordingly by 1.5.3 for any $y\in X$
there exist positive integers $s$ and $r$ with
$r > \frac{s+n}{\eps+\delta}$ such that $|K_X+rL|$ does not separate
$s-$jets at  $y$. Clearly we can assume
$2\eps\ge \delta +\eps > \frac{s+n}{r}\ge\frac{\eps}{2}$, and that $\delta$
is so small that
$$\gamma\ge \gamma':=\frac{\eps\gamma}{\delta+\eps}\ge
\frac{\gamma-1}{2}+1.$$
Since by assumption
$$(rL)^n \ge r^n(\eps\gamma)^n >(s+n)^n\cdot
(\gamma')^n,$$
we can apply Theorem 3.1  to $\gamma'$,
and $\alpha'_i=
\alpha_i\frac{\eps(\gamma-1)-\delta}{(\delta+\eps)(\gamma-1)}$.
By construction
$$2\eps\gamma\ge \gamma'\frac{s+n}{r}\ge\frac{\eps}{2}
\left(\frac{\gamma-1}{2}+1\right),$$
so that from Corollary 3.2 we deduce the existence of a constant $\kappa>0$
and a proper subvariety $V\subset X$ containing $y$
of degree
$$\text{deg}_LV=L^d.V < \left(\frac{n+s}{r}\right)^d
\frac{(\gamma')^n}{(\alpha_d')^{n-d}}\left(\frac{1}{1-\kappa}\right)$$
independently of the choice of $\delta$.
For $\delta>0$ small enough this
gives the existence of
a subvariety $V$ containing $y$ such that
$$L^d.V < \frac{\eps^d\cdot \gamma^n}{\alpha_d^{n-d}},$$
leading to a contradiction.
\qed\enddemo

\noindent
{}From the above Theorem one can deduce easily various
boundedness statements by specifying the $\alpha_i$.
\remark{Example 3.5}  Let $X$ be a smooth projective
variety of dimension $n$, $L$ a nef and big line bundle on $X$
and $x\in X$ a general point.

\noindent
$a)$ Setting $\alpha_1=\dots=\alpha_{n}=1$,
one obtains that the Seshadri constant
satisfies
the universal bound
$$\eps(L,x)\ge \frac{1}{n^n}.$$

\noindent
$b)$  With  the following choice one comes closer to the bound
obtained in \cite{EKL}: put
$$\alpha_i=\frac{n-1}{2^i(1-2^{1-n})},\;\;
\text{ and define }$$
$$\mu(d):= \min_{V}\left\{L^d\cdot V\right\},$$
where the minimum runs over all $d$-dimensional subvarieties
$V\subseteq X$ containing very general points. Then
$$\eps(L,x)\ge\min_{1\le d\le n}
\left\{{\root {\scriptstyle d}\of{\mu(d)
\frac{(n-1)^{n-d}}{n^n[2^d(1-2^{1-n})]^{n-d}}}}\right\}.$$
\endremark
\vglue20pt
\head
\S4 Bounds for Seshadri constants at arbitrary points
\endhead
In this section, we show how to apply the strategy of \S3 to
obtain certain bounds for Seshadri constants
at arbitrary points using the following result
of \cite{ELN}:

\proclaim{Theorem 4.1}(Ein--Lazarsfeld--Nakamaye)
Let $X$ be a smooth $n$--dimensional variety with tangent bundle $T_X$,
$x\in X$ any point,
$A$ an ample line bundle on $X$ and $\delta\ge 0$ a real number
such that $T_X(\delta A)$ ist nef.
Let $$0=\beta_1 <\dots < \beta_{n+1}< {\root n\of{L^n}}$$ be rational
numbers.
Then either

\noindent
$(a)\;\;\;$ there exists $E\in |kA|\;\;(k\gg 0)$ with an almost
isolated singularity at $x$
of index $$\frac{\beta_{n+1}-\beta_n}{1+\delta\beta_n}, \;\;\text{or}$$

\noindent
$(b)\;\;\;$ there exists an irreducible
subvariety $V\ni x$ of codimension $c\ne n$ with
$$A^{n-c}V=\text{deg}_AV\le\left(\frac{1+\delta\beta_c}{\beta_{c+1}-
\beta_c}\right)^c\cdot \beta_{n+1}^n.$$
\endproclaim

Theorem 4.1 follows from \cite{ELN, Theorem 3.9} together with
\cite{ibid, 1.5, 1.6}
and the remark (from the proof of Theorem 3.9) that
in any event $V$ is an irreducible component of the base locus
of the linear system $|\I_V^{(k\eps)}\otimes \O_X(k(1+\delta\sigma)A))|$
on $X$,
which gives, in case $V$ is $0-$dimensional, by Bertini's Theorem
the divisor in (b).

Then the argument proceeds as in \S3:
\proclaim{Corollary 4.2}
Let $X, A$ and $\delta$ be as in (4.1), moreover $r,s$ be positive
integers, and $\gamma >1$ rational such that
$(rA)^n >(\gamma(n+s))^n$. Let
 $\alpha_1,\dots,\alpha_{n-1}$ be positive rational numbers
with $$\left(1+\left(1+\frac{\delta}{r}(n+s)\right)\cdot
\sum_{i=1}^{n-1}\alpha_i\right)=\gamma.$$

Then for any $x\in X$ either $K_X+rA$ separates
$s-$jets at $x$, or there exists a subvariety $V\subseteq X$ through $x$
of dimension $d\ge 1$ and degree
$$A^{n-c}\cdot V\le
\frac{\gamma^n(n+s)^d}{r^d}
\left(\frac{1+\frac{\delta}{r}
(n+s)\sum_{i=d+1}^{n-1}\alpha_i}{\alpha_d}\right)^{n-d}.$$
\endproclaim
\demo{Proof} Put $A'=rA, \beta_{n+1}=\gamma(n+s), \beta_n=(n+s)\cdot
\sum_{i=1}^{n-1}\alpha_i$, and $\beta_i=\beta_{i+1}-\alpha_{n-i}(n+s)$.
Then apply the Theorem 4.1 to $A'$ and $\delta'=\frac{\delta}{r}$,
and use (1.5.1).
\qed\enddemo
\proclaim{Theorem 4.3}
Let $X,A$ and $\delta$ be as in (4.1), moreover $\eps >0$
real and $\alpha_1,\dots, \alpha_n$ be positive rational
numbers.
Let $x\in X$ be any point and suppose that any $d-$dimensional
$(1\le d\le n)$ subvariety $V\subseteq X$ containing $x$
satisfies
$$\text{deg}_AV=A^d\cdot V\ge \eps^d\left(\frac{1+\delta\eps
\sum_{i=d+1}^{n-1}\alpha_i}{\alpha_d}\right)^{n-d}\cdot
\left(1+(1+\delta\eps)\sum_{i=1}^{n-1}\alpha_i\right)^n.$$
Then $\eps(A,x)\ge \eps$.
\endproclaim
\demo{Proof}
Fix $x\in X$ and suppose $\eps(A,x)< \eps$. Then there exist positive
integers $r$ and $s$ with $r >\frac{s+n}{\eps}$ such that $|K_X+rA|$ does
not generate $s-$jets at $x$.
Put $\gamma= 1+\left(1+\frac{\delta}{r}(n+s)\right)\sum_{i=1}^{n-1}\alpha_i$.
Then by assumption
$$\aligned
(rA)^n &\ge r^n\eps^n\left(1+(1+\delta\eps)\sum_{i=1}^{n-1}\alpha_i\right)^n\cr
&> (s+n)^n\left(1+\left(1+
\frac{\delta}{r}(s+n)\right)\sum_{i=1}^{n-1}\alpha_i\right)^n\cr
&= \left((s+n)\gamma\right)^n.\endaligned$$
Therefore Corollary 4.2 gives the existence of $V\ni x$ of dimension $d\ge 1$
with
$$\aligned
A^d\cdot V&\le
\frac{\gamma^n(n+s)^d}{r^d}
\left(\frac{1+\frac{\delta}{r}(n+s)
\sum_{i=d+1}^{n-1}\alpha_i}{\alpha_d}\right)^{n-d}
\cr
&<
\eps^d\gamma^n
\left(\frac{1+\delta\eps\sum_{i=d+1}^{n-1}\alpha_i}{\alpha_d}\right)^{n-d}
\cr
&< \eps^d\left(\frac{1+\delta\eps
\sum_{i=d+1}^{n-1}\alpha_i}{\alpha_d}\right)^{n-d}\cdot
\left(1+(1+\delta\eps)\sum_{i=1}^{n-1}\alpha_i\right)^n.
\endaligned$$
contradicting the assumptions.
\qed\enddemo
Setting $\alpha_1=\dots=\alpha_n=1$  one then obtains:
\proclaim{Corollary 4.4}
Let $X, A$ and $\delta > 0$ be as in (4.1), and
$x\in X$ be any point.
Then
$$\eps(A,x)\ge \text{min}\left\{\frac{1}{(n-1)^{n-1}(2n-1)^n},
\frac{1}{\delta}\right\}.\hfil\qed$$
\endproclaim
\remark{Remark 4.5}
It is well known that, for a very  ample line bundle
$H$ on $X$, the twisted tangent bundle
$T_X(K_X+nH)$ is globally generated, and in particular nef.
In case $A=K_X$ is ample on $X$ one therefore can use one of
the available effectivity
statements for very ampleness of ample line bundles
(e.g. \cite{De93}) to determine explicit values for $\delta$,
making Theorem 4.3 or Corollary 4.4
effective in the way that the bounds for the Seshadri constant
of $K_X$ at any $x\in X$ do
only depend on the dimension $n$.
The same argument works in case $A=-K_X$ is ample, or
for any ample $A$ in case $K_X$ is trivial.
\endremark
\remark{Remark 4.6}
Let us finally  compare Corollary 4.4 with the bounds
that can be obtained using Angehrn--Siu's basepoint--free
Theorem. Namely,
Angehrn and Siu prove that, for an ample line bundle $A$
on $X$, the adjoint line bundles $mA+K_X$ are free for
$m\ge \frac{1}{2}n(n+1) +1$. An elementary argument
(cf. e.g. \cite{K\"u, 3.3}) shows
that $\eps(A,x)\ge 1$ for all ample basepoint--free  line bundles
$A$. Moreover, if $T_X(\delta A)$ is nef,
so is the ${\Bbb Q}-$line bundle $M:=\text{det}\left(T_X(\delta A)\right)
=-K_X+n\delta A$.

By definition  Seshadri constants have a sublinearity property saying
that, for any rational $\lambda, \mu\ge 0$ and nef line bundles $L$ and $M$,
the inequality $\eps(\lambda L+\mu M,x)\ge \lambda\cdot\eps(L,x)
+\mu\cdot\eps(M,x)$
holds.  This shows
$$\eps(A,x)\ge \frac{2}{n(n+2\delta+1)+2}.$$
\endremark
\vfill

\enddocument
\end